\documentclass{amsart}
\usepackage{mathrsfs}
\usepackage{graphicx}
\usepackage{amsmath}
\usepackage{amscd}
\usepackage{amssymb}

\input xy
\xyoption{all}

\newcommand{\pd}{\partial}
\newcommand{\bC}{{\mathbb C}}

\newcommand{\bR}{{\mathbb R}}
\newcommand{\bZ}{{\mathbb Z}}

\newcommand{\half}{\frac{1}{2}}

\newtheorem{theorem}{Theorem}[section]
\newtheorem{theorem/definition}{Theorem/Definition}[section]

\newtheorem{prop}{Proposition}[section]

\theoremstyle{remark}

\theoremstyle{definition}

\newcommand{\be}{\begin{equation}}
\newcommand{\ee}{\end{equation}}
\newcommand{\bea}{\begin{eqnarray}}
\newcommand{\ben}{\begin{eqnarray*}}
\newcommand{\een}{\end{eqnarray*}}
\newcommand{\eea}{\end{eqnarray}}
\newcommand{\bet}{\begin{equation}
\begin{split}}
\newcommand{\eet}{\end{split}
\end{equation}}

\usepackage{color}

\definecolor{yellow}{rgb}{1,1,0}
\definecolor{orange}{rgb}{1,.7,0}
\definecolor{red}{rgb}{1,0,0} \definecolor{green}{rgb}{0,1,1}
\definecolor{white}{rgb}{1,1,1}

\definecolor{A}{rgb}{.75,1,.75}

\begin{document}

\title
{Hessian Geometry and Phase Changes of Multi-Taub-NUT Metrics}

\author{Jian Zhou}
\address{Department of Mathematical Sciences\\
Tsinghua University\\Beijing, 100084, China}
\email{jzhou@math.tsinghua.edu.cn}

\begin{abstract}
We study the Hessian geometry of toric multi-Taub-NUT metrics  and their phase change phenomena
via the images of their moment maps.
This generalizes an earlier paper on toric Gibbons-Hawking metrics.
\end{abstract}
\maketitle

%% \tableofcontents

\section{Introduction}

This is a sequel to \cite{Zhou-Kepler} and \cite{Zhou-GH}
in which the Kepler metrics and the toric Gibbons-Hawking metrics
are studied from the point of view of Hessian geometry \cite{Shima}, respectively.
These are  techniques developed in the setting of toric canonical K\"ahler metrics
on compact toric manifolds in K\"ahler geometry \cite{Gui, Abr, Don},
and also in the the study of toric Sasaki-Einstein metrics \cite{GMSW, Mar-Spa, Mar-Spa-Yau}
that arise in AdS/CFT correspondence in string theory \cite{Mal, Wit}.
Our motivation is to test the applicability of some techniques developed in
string theory to problems in classical gravity and general relativity,
and to gain some new insights to these problems in doing this.
For the examples in our earlier papers mentioned above,
by examining the images of the moment maps,
one can easily visualize the degenerations of K\"ahler metrics,
more importantly, 
one can study their phase changes \cite{Zhou-GH, Wang-Zhou},
a concept borrowed from statistical physics and 
introduced in \cite{Duan-Zhou1, Duan-Zhou2} in the context of K\"ahler geometry.
In this paper we will show that the results in \cite{Zhou-GH} can be generalized
to toric multi-Taub-NUT spaces with mild modifications,
and more importantly,
we will report a new type of phase transition not known
from our earlier papers \cite{Duan-Zhou1, Duan-Zhou2, Zhou-GH, Wang-Zhou}.

Both the Gibbons-Hawking metrics and the multi-Taub-NUT metrics are
gravitational instantons in Euclidean gravity in dimension four of type $A_{n-1}$.
They share the same construction that shows that they
are both hyperk\"ahler with a triholomorphic circle action.
They live on the same space: The crepant resolution of $\widehat{\bC^2/\bZ_n}$,
but with different asymptotic behavior near the infinity:
The Gibbons-Hawking metrics are asymptotically locally Euclidean (ALE)
while the multi-Taub-NUT metrics are asymptotically locally flat (ALF).
This means they have the following asymptotic form respectively:
\bea
&& g_{n-GH} \sim dr^2 + r^2(\sigma_1^2+ \sigma_2^2+\sigma_3^2), \\
&& g_{n-TB} \sim dr^2 + \sigma_1^2 + r^2(\sigma_2^2+\sigma_3^2),
\eea
where $\sigma_1, \sigma_2, \sigma_3$ are left-invariant one-forms on $S^3/\bZ_n$.

As in the case of Gibbons-Hawking metrics,
we will focus on the toric  multi-Taub-NUT spaces
and consider the convex bodies that arise as images of natural torus actions on them.
Our main result is that the metrics and the complex structures can be
reconstructed from some functions on the convex bodies.
Considerations of the degenerations of these convex bodies
lead to the consideration of phase transition phenomena
of the multi-Taub-NUT metrics.
We will also consider another kind of phase transition.
Recall the Gibbons-Hawking metrics and the multi-Taub-NUT metrics
are given by the same construction,
the only difference is a positive parameter in the latter.
When this parameter becomes zero,
the multi-Taub-NUT metrics become the Gibbons-Hawking metrics.
We will let the parameter to further become negative.
Then the total space is divided into two regions,
in one of which the metrics remain positive definite,
but on the other, the metrics become {\em negative} definite.

The rest of the paper is arranged as follows.
In \S \ref{sec:Moment},
after recalling the Gibbons-Hawking construction  of the multi-Taub-NUT metrics,
we study the moment maps of the toric cases and the phase transition of their images.
Hessian geometry is used in \S \ref{sec:Hessian} to explicitly
construct local complex coordinates on multi-Taub-NUT spaces.
We describe the phase change for $1$-Taub-NUT metrics in \S \ref{sec:n=1}.
We end the paper by some concluding remarks in \S \ref{sec:Conclusions}.

\section{Phase Transitions of Moment Map Images of Toric Multi-Taub-NUT Metrics}
\label{sec:Moment}

In this Section we first recall the Gibbons-Hawking construction \cite{Gib-Haw}
of multi-Taub-NUT metrics,
then we focus on the toric cases and consider their moment maps.
We explain
how the moment map images undergo a phase transition.

\subsection{Gibbons-Hawking construction of multi-Taub-NUT metrics}

Given $n$ distinct points $\vec{p}_1, \dots, \vec{p}_n$ in $\bR^3$,
consider a function $V$ defined by
\be
V_\epsilon(\vec{r}) =  \frac{\epsilon}{2} + \half \sum_{j=1}^n \frac{1}{|\vec{r}-\vec{p}_j|},
\ee
for some fixed constant $\epsilon \in \bR$.
Since  $V_\epsilon$ is a  harmonic function
on a simply connected domain $U\subset \bR^3-\{\vec{p}_1, \dots, \vec{p}_n\}$ be a domain
one has
\be
*dV_\epsilon = - d\alpha
\ee
for some smooth one-form $\alpha$,
where $*$ is the Hodge star-operator.
As in the case of  toric Gibbons-Hawking metrics \cite[\S 2.2]{Zhou-GH},
we will let the points $\vec{p}_1, \dots, \vec{p}_n$ lie in a line, say,
$\vec{p}_j = (0, 0, c_j)$, $j=1, \dots, n$,
$c_1 < c_2 < \cdots < c_n$,
we take $U:=\bR^3-\{(0,0,z)\;|\;z \geq c_1\}$ we take:
\be \label{eqn:alpha}
\alpha
= - \half \sum_{j=1}^n \frac{xdy-ydx}{r_j(r_j-z+c_j)},
\ee
where $r_j = \sqrt{x^2+y^2+(z-c_j)^2}$.

On the principal bundle $U \times S^1 \to U$ with connection $1$-form $d\varphi+\alpha$,
where $\varphi$ is the natural coordinate on $S^1$, i.e., $e^{i\varphi} \in S^1$,
consider
\be
g_\epsilon = \frac{1}{V_\epsilon}(d\varphi + \alpha)^2 + V_\epsilon \cdot (dx^2+dy^2+dz^2).
\ee
For $\epsilon > 0$,
this is
the {\em multi-Taub-NUT metric} in local coordinates $\{\varphi, x, y, z\}$;
for $\epsilon = 0$,
this is the {\em Gibbons-Hawking metric};
for $\epsilon < 0$,
$g_\epsilon$ defines a Riemannian metric in  the region $\Omega_+$ where $V_\epsilon > 0$,
and $-g_\epsilon$ defines a Riemannian metric in the region $\Omega_-$ where $V_\epsilon < 0$.

\subsection{The regions $\Omega_+$ and $\Omega_-$ for $\epsilon < 0$}

Let $U_\pm$ be the regions in $\bR^3-\{\vec{p}_1, \dots, \vec{p}_n\}$
on which $V_\epsilon$ takes the $\pm$-sign,
so that $\Omega_\pm$ are circle bundles over them.

Write $\epsilon = - a$ for some $a > 0$.
Let us first consider the case of $n=1$ and take $c_1 = 0$.
Then we have
\ben
&& V_\epsilon = - \frac{a}{2}+ \frac{1}{2\sqrt{\rho^2+z^2}},
\een
where $\rho = \sqrt{x^2+y^2}$.
When $V_\epsilon >0$, $0 \leq \rho^2 < \frac{1}{a^2} -z^2$.
When $V_\epsilon < 0$, $\rho^2 > \frac{1}{a^2} -z^2$.
This imposes no bounds on $z$.
\emph{}

Now we move on to  the general case. Fix $z$,
and regard $V_{-a}$ as function of $p= \rho^2$,
\ben
&& V_{-a} = - \frac{a}{2}+ \half \sum_{j=1}^n \frac{1}{\sqrt{p+(z-c_j)^2}},
\een
then we have
\be
\frac{\pd V_{-a}}{\pd p} = - \frac{1}{4} \sum_{j=1}^n \frac{1}{\sqrt{(p+(z-c_j)^2)^3}},
\ee
Therefore, there exists a smooth function $p_a(z)$ such that
when $\rho = \sqrt{p_a(z)}$,
$V_{-a} = 0$,
when $\rho < \sqrt{p_a(z)}$, $V_{-a} > 0$,
and when $\rho > \sqrt{p_a(z)}$, $V_{-a} < 0$.
In other words,
$U_+$ is given by $\rho^2 < p_a(z)$
and $U_-$ is given by $\rho^2 > p_a(z)$.

\subsection{Complex structures}
Consider the almost complex structure given by:
\begin{align}
J_\epsilon^*(V_\epsilon^{-1/2}(d\varphi + \alpha)) & = - V_\epsilon^{1/2}dz, &
J_\epsilon^*(V_\epsilon^{1/2} d z) & = V_\epsilon^{-1/2}(d\varphi + \alpha), \\
J_\epsilon^*(V_\epsilon^{1/2}dx) & = - V_\epsilon^{1/2} dy, &
J_\epsilon^*(V_\epsilon^{1/2}dy) & = V_\epsilon^{1/2} dx.
\end{align}
As in \cite{Zhou-GH},
one can check that it is integrable,
with the space of type $(1, 0)$-forms  generated by:
\begin{align}
dx + \sqrt{-1} dy, && (d\varphi + \alpha) + \sqrt{-1}V_\epsilon dz.
\end{align}

When $\epsilon < 0$,
$J_\epsilon$ is undefined when $V_\epsilon = 0$,
and one has to consider $J_\epsilon$ in $\Omega_+$
and $\Omega_-$ separately.

\subsection{Symplectic structures}

The complex structure $J_\epsilon$ is compatible with the Riemannian metric $g_\epsilon$,
with the symplectic form given by:
\be \label{eqn:Symplectic}
\omega_\epsilon = (d\varphi+\alpha) \wedge d z + V_\epsilon d x \wedge dy.
\ee
Of course when $\epsilon < 0$,
$\omega_\epsilon$ is degenerate along $V_\epsilon = 0$.
As in the case of Gibbons-Hawking metric (i.e., $g_0$),
the multi-Taub-NUT metric (i.e. $g_\epsilon$ for $\epsilon>0$) is hyperk\"ahler.
For $\epsilon< 0$,
$g_\epsilon$ is hyperk\"ahler in the region $\Omega_+$ where $V_\epsilon > 0$,
and $-g_\epsilon$ is hyperk\"ahler in the region $\Omega_-$ where $V_\epsilon < 0$.

\subsection{Torus action and moment map on a toric multi-Taub-NUT space}

Using the explicit choice \eqref{eqn:alpha} of $\alpha$ in $U$ we have
\be \label{eqn:Symplectic2}
\omega_\epsilon = \biggl(d\varphi- \half \sum_{j=1}^n \frac{xdy-ydx}{r_j(r_j-z+c_j)}\biggr) \wedge d z
+ \half\biggl(\epsilon + \sum_{j=1}^n \frac{1}{r_j} \biggr) d x \wedge dy.
\ee
Let $x+y\sqrt{-1} y = \rho e^{\sqrt{-1}\theta}$,
then
\ben
\omega_\epsilon
& = & \biggl(d\varphi- \half \frac{\rho^2d\theta}{r_j(r_j-z+c_j)}\biggr) \wedge d z
+ \half\biggl(\epsilon + \sum_{j=1}^n \frac{1}{r_j} \biggr) \rho d \rho \wedge d\theta \\
& = & -dz \wedge d\varphi
+ \half \biggl( \frac{r_j+z-c_j}{r_j}   d z
+ \biggl( \epsilon + \sum_{j=1}^n \frac{1}{r_j} \biggr) \rho d \rho \biggr) \wedge d\theta \\
& = & d\mu_1 \wedge d \theta_1 + d\mu_2 \wedge d\theta_2,
\een
where
\begin{align}
\theta_1 & = \varphi, & \mu_1 & = - z, \\
\theta_2 & = \theta, & \mu_2 & = \frac{1}{4} \epsilon (x^2+y^2) + \half\sum_{j=1}^n (r_j+z-c_j).
\end{align}
Here we have used:
\bea
d \mu_2 %%%% & = & \half\sum_{i=1}^n (dr_i + dz)    \\
& = & \half \biggl(\epsilon(xdx+ydy) + \sum_{j=1}^n (\frac{xdx+ydy+(z-c_j)dz}{r_j} + dz)\biggr). \label{eqn:dmu2}
\eea
The Hamiltonian functions $\mu_1$ and $\mu_2$ generate a $2$-torus
action given in local coordinates $(\varphi, x, y, z)$ by:
\be \label{eqn:action}
(e^{i\theta_1}, e^{i\theta_2})  \cdot (\varphi, x, y, z)
= (\varphi+ \theta_1, x\cos\theta_2 - y \sin \theta_2, x\sin \theta_2+y\cos \theta_2, z).
\ee

Note for $\epsilon < 0$,
even though $\omega_\epsilon$ is not defined everywhere,
the group action and the moment map are defined everywhere.
The group action is smooth everywhere, but the moment map
is not.

\subsection{Phase transition of images of the moment maps}

In this subsection,
we explain how the images of the moment maps undergo a phase change when
the parameter $\epsilon$ changes from a positive number to a negative number.

When $\epsilon \geq 0$,
note
\ben
\mu_2& = & \frac{\epsilon}{4}(x^2+y^2) + \half \sum_{j=1}^n (\sqrt{x^2+y^2+(z-c_j)^2} + z-c_j) \\
& \geq & \half \sum_{j=1}^n (|z-c_j|+z-c_j).
\een
Since
\be
|z-c_j|+z-c_j = \begin{cases}
0, & \text{if $z\leq c_j$}, \\
2(z-c_j), & \text{if $z \geq c_j$},
\end{cases}
\ee
it is easy to see that the image of the moment map is the convex region given by
the following inequalities:
\ben
&& l_0: = \mu_2 \geq 0, \\
&& l_1: = \mu_2 + (\mu_1 + c_1) \geq 0, \\
&& l_2: = \mu_2 + (\mu_1 + c_1) + (\mu_1 + c_2) \geq 0, \\
&&  \cdots\cdots \cdots \cdots \cdots \cdots \\
&& l_n: = \mu_2 + \sum_{j=1}^n (\mu_1 + c_j) \geq 0.
\een
This is the same as in the case of toric Gibbons-Hawking metrics.

Now we consider the case of $\epsilon < 0$.
Let us first consider the case of $n=1$ and take $c_1 = 0$.
Write $\epsilon = - a$ for some $a > 0$.
Then we have
\ben
&& \mu_2 = -\frac{a}{4}\rho^2 + \half (\sqrt{\rho^2+z^2} + z).
\een
Consider the function
$$f_z(p) = - \frac{a}{4}p + \half (\sqrt{p+z^2}+z),$$
as a function of $p$.
We have
\ben
f'_z(p) = - \frac{a}{4} + \frac{1}{4} \frac{1}{\sqrt{p+z^2}}
= \half V_{-a}.
\een
So we have $f'_z(x) > 0$ for $0 \leq x < \frac{1}{a^2}-z^2$, and $f_z(0) = \half(|z|+z)$,
$f_z(\frac{1}{a^2}-z^2) = \frac{1}{4a} (1+az)^2$,
it follows that in $\Omega_+$ we have:
\be
\half (|z|+z) \leq \mu_2 < \frac{1}{4a} (1+az)^2.
\ee
Hence moment map image  of $\Omega_+$ is given by the following inequalities:
\be
\half (|\mu_1|-\mu_1) \leq \mu_2 < \frac{1}{4a} (1-a\mu_1)^2.
\ee
This is no longer a convex set.
Its boundary has three pieces:
an interval on the line $\mu_2 = 0$,
an interval on the line $\mu_2 = -\mu_1$,
and they are both tangent to a portion of the curve $\mu_2 = \frac{1}{4a}(1-a\mu_1)^2$
for $-\frac{1}{a} < \mu_1 < \frac{1}{a}$.

We have $f'_z(x) < 0$ for $x \geq \frac{1}{a^2}-z^2$, and
$f_z(\frac{1}{a^2}-z^2) = \frac{1}{4a} (1+az)^2$,
so we have:
\be
\mu_2 < \frac{1}{4a} (1+az)^2.
\ee
The moment image of the region $\Omega_-$ is given by:
\be
\mu_2 < \frac{1}{4a} (1-a\mu_1)^2.
\ee

In this case we obtain an unexpected convexity result as follows:
The union of the moment image of $\Omega_+$ with the complement of the moment
image of $\Omega_-$ is convex.

Now we move on to  the general case:
\ben
&& \mu_2 = -\frac{a}{4}\rho^2 + \half \sum_{j=1}^n (\sqrt{\rho^2+(z-c_j)^2} + z-c_j).
\een
Consider the function
$$f_z(p) = - \frac{a}{4}p + \half \sum_{j=1}^n (\sqrt{p+(z-c_j)^2}+z-c_j),$$
as a function of $p$.
We have
\ben
&& f'_z(p) = - \frac{a}{4} + \frac{1}{4}\sum_{j=1}^n \frac{1}{\sqrt{p+(z-c_j)^2}}= \half V_{-a}, \\
&& f''_z(p) = -\frac{1}{8}\sum_{j=1}^n \frac{1}{\sqrt{(p+(z-c_j)^2)^3}} < 0.
\een
Recall that, when $0 \leq \rho=\sqrt{p_a(z)}$,
$V_{-a} = 0$,
when $\rho < \sqrt{p_a(z)}$, $V_{-a} > 0$,
and so $\mu_2$ is bounded between $f_z(0) = \half \sum_{j=1}^n (|z-c_j|+z-c_j)$
and $f_z(p_a(z)) =
-\frac{a}{4}p_a(z) + \half \sum_{j=1}^n (\sqrt{p_a(z)+(z-c_j)^2} + z-c_j)$.
Since $z= - \mu_1$,
the moment image of $\Omega_+$
is given by:
\be
\begin{split}
&\half \sum_{j=1}^n (|\mu_1+c_j|-\mu_1-c_j) \leq
\mu_2 \\
& < -\frac{a}{4}p_a(-\mu_1)
+ \half \sum_{j=1}^n (\sqrt{p_a(-\mu_1)+(-\mu_1-c_j)^2} -\mu_1 -c_j).
\end{split}
\ee
When $\rho > \sqrt{p_a(z)}$, $V_{-a} < 0$,
$\mu_2 < f_z(p_a(z))$,
and the moment image of the region $\Omega_-$ is
\be
\mu_2 < -\frac{a}{4}p_a(-\mu_1) + \half \sum_{j=1}^n (\sqrt{p_a(-\mu_1)+(-\mu_1-c_j)^2} -\mu_1-c_j).
\ee
We conjecture that the boundary of this region is a convex curve which is concave up
and tangent to the boundary of the moment image in the $\epsilon > 0$ case.

\section{Hessian Geometry of Toric Multi-Taub-NUT Spaces}
\label{sec:Hessian}

In this Section we find the complex potential functions
of the toric multi-Taub-NUT metrics and use them to define local complex coordinates.

\subsection{Symplectic coordinates and Hessian geometry for toric
multi-Taub-NUT spaces}

The following result can be obtained by a simple modification of the corresponding
result in \cite{Zhou-GH}:

\begin{prop}
In the above symplectic coordinates,
the multi-Taub-NUT metric takes the following form:
\be
g_\epsilon = \sum_{i,j=1}^2 (\frac{1}{2} G_{ij} d\mu_i d\mu_j +2 G^{ij}d\theta_id\theta_j),
\ee
where the coefficient matrices $(G_{ij})_{i,j=1, 2}$ and $(G^{ij})_{i,j=1,2}$ are given by:
\ben
(G_{ij})_{i,j=1, 2}
= \begin{pmatrix}
2V_\epsilon  +
\frac{\rho^2}{2V_\epsilon} \biggl(\sum_{j=1}^n\frac{1}{r_j(r_j-(z-c_j))}\biggr)^2 &
\frac{1}{V_\epsilon} \sum_{j=1}^n\frac{1}{r_j(r_j-(z-c_j))} \\
\frac{1}{V_\epsilon} \sum_{j=1}^n\frac{1}{r_j(r_j-(z-c_j))} &
\frac{2}{V_\epsilon\rho^2}
\end{pmatrix}
\een

\ben
( G^{ij})_{i,j=1,2}=
\begin{pmatrix}
\frac{1}{2V_\epsilon} & -\frac{1}{4V_\epsilon} \sum_{j=1}^n\frac{\rho^2}{r_j(r_j-(z-c_j))}  \\
-\frac{1}{4V_\epsilon} \sum_{j=1}^n\frac{\rho^2}{r_j(r_j-(z-c_j))}
& \frac{V_\epsilon \rho^2}{2} +
\frac{1}{8V_\epsilon} \biggl(\sum_{j=1}^n\frac{\rho^2}{r_j(r_j-(z-c_j))}\biggr)^2
\end{pmatrix}
\een
These matrices are inverse to each other.
The complex potential and the K\"ahler potential are given by the following formulas respectively:
\be \label{eqn:Complex}
\begin{split}
\psi_\epsilon & = \half \sum_{j=1}^n \biggl( (r_j + (z-c_j)) \log (r_j + (z-c_j)) \\
& + (r_j-(z-c_j)) \log(r_j-(z-c_j)) \biggr)
+ \frac{\epsilon}{2} z^2 +C_1\mu_1 + C_2\mu_2.
\end{split}
\ee
And the K\"ahler potential is given by:
\be \label{eqn:Kahler}
\psi_\epsilon^\vee= -\sum_{j=1}^n c_j \log(r_j-(z-c_j))
+ \frac{\epsilon}{2} z^2 +C_1\mu_1 + C_2\mu_2.
\ee
for some constants $C_1$ and $C_2$.
\end{prop}

\subsection{Hessian local complex coordinates}

The function $\psi$ is called the {\em complex potential}
because one can find local complex coordinates $z_1$ and $z_2$ so that
\be
\frac{dz_i}{z_i} = \half \sum_{j=1}^2 \frac{\pd^2\psi_\epsilon}{\pd \mu_i \pd \mu_j} d\mu_j
+ \sqrt{-1} d\theta_i =  \half \sum_{i,j=1}^2 G_{ij} d\mu_j + \sqrt{-1} d \theta_i
\ee
is of type $(1, 0)$.
We have
\bea
\frac{dz_1}{z_1} & = &
- \frac{\epsilon}{2} dz  + \half d \sum_{j=1}^n  \log (r_j - (z-c_j))
+ \sqrt{-1} d\theta_1, \label{eqn:dz1/z1} \\
\frac{dz_2}{z_2} & = & d \log \rho + \sqrt{-1} d\theta_2 = d \log (x + y \sqrt{-1}).
\label{eqn:dz2/z2}
\eea
Therefore, we take
\bea
&& z_1 = \prod_{j=1}^n (r_j-(z-c_j))^{1/2}\cdot e^{- \frac{\epsilon}{2} z +\sqrt{-1} \theta_1},
\label{eqn:z1} \\
&& z_2 = x+ \sqrt{-1} y. \label{eqn:z2}
\eea

A simple modification by changing $V$ to $V_\epsilon$
in the proof of Theorem 4.1 in \cite{Zhou-GH} then proves the following:

\begin{theorem}
The metrics $g_\epsilon$ and K\"ahler forms $\omega_\epsilon$ are given in local complex coordinates $z_1, z_2$
as follows:
\be \label{eqn:g-in-complex}
\begin{split}
g_\epsilon = &   \frac{1}{V_\epsilon} \frac{dz_1}{z_1}\frac{d\bar{z}_1}{\bar{z}_1}
- \frac{1}{2V_\epsilon} \sum_{j=1}^n\frac{r_j+z-c_j}{r_j}
\biggl( \frac{dz_1}{z_1} \frac{d\bar{z}_2}{\bar{z}_2}
+ \frac{dz_2}{z_2} \frac{d\bar{z}_1}{\bar{z}_1} \biggr) \\
+ & \biggl[ V_\epsilon \rho^2 + \frac{1}{4V}
\biggl( \sum_{j=1}^n \frac{r_j+z-c_j}{r_j}\biggr)^2 \biggr]\emph{}
\frac{dz_2}{z_2} \frac{d\bar{z}_2}{\bar{z}_2} ,
\end{split}
\ee
\be
\begin{split}
\omega_\epsilon & = \frac{1}{2\sqrt{-1}} \biggl(\frac{1}{V_\epsilon} \frac{dz_1}{z_1} \wedge
\frac{d\bar{z}_1}{\bar{z}_1} \\
& - \frac{1}{2V_\epsilon} \sum_{j=1}^n\frac{r_j+z-c_j}{r_j}
\biggl( \frac{dz_1}{z_1} \wedge \frac{d\bar{z}_2}{\bar{z}_2}
+ \frac{dz_2}{z_2}\wedge  \frac{d\bar{z}_1}{\bar{z}_1} \biggr) \\
& + \biggl[ V_\epsilon \rho^2 + \frac{1}{4V_\epsilon}
\biggl( \sum_{j=1}^n \frac{r_j+z-c_j}{r_j}\biggr)^2 \biggr]
\frac{dz_2}{z_2} \wedge \frac{d\bar{z}_2}{\bar{z}_2} \biggr).
\end{split}
\ee
\end{theorem}

\subsection{The $(\alpha,\beta)$-coordinates}

As in \cite{Zhou-GH},
we make the following change of variables:
\be \label{eqn:Alpha-Beta}
z_1 = \beta_1,  \;\;\;\;\; z_2 = \alpha_1\beta_1.
\ee
By a simple calculation we have:

\begin{theorem}
The metrics $g_\epsilon$  are given in local complex coordinates $\alpha_1, \beta_1$
as follows:
\be \label{eqn:g-in-ab}
\begin{split}
g_\epsilon = & \biggl[ V_\epsilon \rho^2 + \frac{1}{4V_\epsilon}
\biggl( \sum_{j=1}^n \frac{r_j+z-c_j}{r_j}\biggr)^2 \biggr]
\frac{d\alpha_1}{\alpha_1}  \cdot \frac{d\bar{\alpha}_1}{\bar{\alpha}_1} \\
+ & \biggl[ V_\epsilon \rho^2 - \frac{1}{2V_\epsilon} \sum_{j=1}^n\frac{r_j+z-c_j}{r_j}
+ \frac{1}{4V_\epsilon} \biggl( \sum_{j=1}^n \frac{r_j+z-c_j}{r_j}\biggr)^2 \biggr] \\
& \cdot \biggl(\frac{d\alpha_1}{\alpha_1}\cdot \frac{d\bar{\beta}_1}{\bar{\beta}_1}
+\frac{d\beta_1}{\beta_1} \cdot \frac{d\bar{\alpha}_1}{\bar{\alpha}_1} \biggr) \\
+ & \biggl[\frac{1}{V_\epsilon} + V_\epsilon \rho^2 - \frac{1}{V_\epsilon} 
\sum_{j=1}^n\frac{r_j+z-c_j}{r_j}
+ \frac{1}{4V_\epsilon} \biggl( \sum_{j=1}^n \frac{r_j+z-c_j}{r_j}\biggr)^2 \biggr]
\frac{d\beta_1}{\beta_1}\cdot \frac{d\bar{\beta}_1}{\bar{\beta}_1}.
\end{split}
\ee
\end{theorem}

As in \cite[\S 4.5]{Zhou-GH},
one can introduce $(\alpha_i, \beta_i)$ for $i=1, \dots, n$, such that
\begin{align} \label{eqn:Coordinate change}
\alpha_{i+1} & = \alpha_i^2 \beta_i, & \beta_{i+1} & = \alpha_i^{-1},
\end{align}
for $i=1, \dots, n-1$. So the underlying space for the family of metrics $g_\epsilon$
is the crepant resolution of $\bC^2/\bZ_n$.

\section{Phase Transition of The $1$-Taub-NUT Metrics with Respect to $\epsilon$}
\label{sec:n=1}

In this Section we discuss the phase transition of the $1$-Taub-NUT metrics
with respect to the parameter $\epsilon$.

\subsection{The complex potential in the $n=1$ case} One can take $c_1=0$.
From the equation
\be
\frac{\epsilon}{4} \rho^2 + \sqrt{\rho^2+z^2} + z = 2 \mu_2
\ee
we can solve for $\rho^2$:
\be
\rho^2 = -4 \frac{-2\epsilon \mu_2-(2-\epsilon z)
+\sqrt{8\epsilon\mu_2+(2-\epsilon z)^2}}{\epsilon^2} = 4\mu_2 (\mu_2-z) + \cdots
\ee
and
\be
r_1 = \frac{-2+\sqrt{8\epsilon\mu_2+(2-\epsilon z)^2}}{\epsilon} = 2\mu_2-z  + \cdots,
\ee
where $\cdots$ are higher order terms in $\epsilon$.

The complex potential is
\be \label{eqn:Complex-n=1}
\begin{split}
\psi & = \frac{\epsilon z-2+\sqrt{8\epsilon\mu_2+(2-\epsilon z)^2}}{2\epsilon}
\log \frac{\epsilon z-2+\sqrt{8\epsilon\mu_2+(2-\epsilon z)^2}}{\epsilon} \\
& + \frac{-\epsilon z-2+\sqrt{8\epsilon\mu_2+(2-\epsilon z)^2}}{2\epsilon}
\log \frac{-\epsilon z-2+\sqrt{8\epsilon\mu_2+(2-\epsilon z)^2}}{\epsilon} \\
& +C_1\mu_1 + C_2\mu_2 + \frac{\epsilon}{2} \mu_1^2,
\end{split}
\ee
for some constants $C_1, C_2$.

\subsection{The $1$-Taub-NUT metrics in $(\alpha, \beta)$-coordinates}

In this case write $\alpha_1=\alpha$ and $\beta_1 = \beta$.
The $1$-Taub-NUT metric becomes:
\be \label{eqn:Metric-n=1}
\begin{split}
g_\epsilon & =  \frac{\epsilon^2r^2+2\epsilon r-\epsilon^2 r z+2-2\epsilon z}{2(1+\epsilon r)}
\cdot e^{-\epsilon z} d\alpha d\bar{\alpha} \\
& + \frac{\epsilon (2+\epsilon r)}{2 (1+\epsilon r)} \cdot \biggl(
\bar{\alpha}\beta d\alpha  d\bar{\beta} + \bar{\beta}\alpha d\beta d\bar{\alpha} \biggr) \\
& +  \frac{\epsilon^2r^2+2\epsilon r+\epsilon^2rz+2+2\epsilon z}{2(1+\epsilon r)}
\cdot e^{\epsilon z} d\beta d\bar{\beta}.
\end{split}
\ee

Let us now analyze the phase transition of the family of metrics \eqref{eqn:Metric-n=1}.
From
\begin{align}
z_1 & = ((x^2+y^2+z^2)^{1/2}-z)^{1/2}\cdot e^{\frac{-\epsilon z}{2}+\sqrt{-1} \theta_1},
& z_2 & = x+ \sqrt{-1} y.
\end{align}
we get:
\bea
&& \rho^2  = |z_2|^2, \\
&& ((|z_2|^2+z^2)^{1/2}-z)e^{-\epsilon z}  = |z_1|^2. \label{eqn:z1-mod}
\eea
Or in the $(\alpha, \beta)$-coordinates
\bea
&& \rho^2 = |\alpha|^2|\beta|^2, \\
&& ((|\alpha|^2|\beta|^2+z^2)^{1/2}-z)e^{-\epsilon z}  = |\beta|^2.
\eea
In order to consider the solution of
the second equation,
let
$$f(z) = ((|\alpha|^2|\beta|^2+z^2)^{1/2}-z)e^{-\epsilon z}$$
be regarded as a function
in $z$ with other variables as parameters,
and consider its derivative:
\ben
f'(z) & = & e^{-\epsilon z}(-\epsilon ((|\alpha|^2|\beta|^2+z^2)^{1/2}-z)
+ \frac{z}{(|\alpha|^2|\beta|^2+z^2)^{1/2}}- 1) \\
& = & -e^{-\epsilon z}((|\alpha|^2|\beta|^2+z^2)^{1/2}-z) (\epsilon
+ \frac{1}{(|\alpha|^2|\beta|^2+z^2)^{1/2}}) \\
& = & -2 V_\epsilon  e^{-\epsilon z}((|\alpha|^2|\beta|^2+z^2)^{1/2}-z).
\een
And so when $\epsilon \geq 0$, $\alpha\beta \neq 0$,
we have $f'(z) < 0$,
and therefore, there is only one $z$ such that \eqref{eqn:z1-mod} holds,
and by inverse function theorem,
$z$ is a smooth function in $|\alpha|^2$ and $|\beta|^2$.
Near $|\alpha|^2=|\beta|^2= \epsilon=0$,
$z$ is analytic in these variables,
and one can see that
\ben
z& = & (|\alpha|^2-|\beta|^2)\biggl(\frac{1}{2}-\frac{|\alpha|^2+|\beta|^2}{4} \epsilon
+\frac{3|\alpha|^4+2|\alpha|^2|\beta|^2+3|\beta|^4}{16} \epsilon^2 + \cdots \biggr), \\
%%% > -(b^2+a^2)*(2*b^4-a^2*b^2+2*a^4)/12*t^3+(125/768*a^8+11/192*a^6*b^2+23/384*a^4*b^4+11/192*a^2*b^6+125/768*b^8)*t^4)
r & = & \frac{|\alpha|^2+|\beta|^2}{2} + (|\alpha|^2-|\beta|^2)^2\biggl(
-\frac{\epsilon}{4} +\frac{3|\alpha|^2+3|\beta|^2}{16} \epsilon^2 + \cdots \biggr)
\een

When $\epsilon < 0$,
write $\epsilon = -a$.
To have $V_\epsilon > 0$,
we need $|z| < \sqrt{\frac{1}{a^2} - |\alpha|^2|\beta|^2}$,
in particular, $|\alpha||\beta| < a$.
So this is the region on the $(\alpha, \beta)$-plane
where \eqref{eqn:Metric-n=1} defines a hyperk\"ahler metric by the same argument as above.
When $V_\epsilon < 0$,
we have $|z| > \sqrt{\frac{1}{a^2} - |\alpha|^2|\beta|^2}$ and $r \geq \frac{1}{a}$,
in particular, $|\alpha||\beta| > a$.
Because $r$ and $z$ are smooth functions in $|\alpha|^2$ and $|\beta|^2$
in this region,
\eqref{eqn:Metric-n=1} defines a negative of a hyperk\"ahler metric there.

To summarize,
we have proved the following:

\begin{theorem}
Suppose that $z$ and $r$ are determined by
\bea
&& ((|\alpha|^2|\beta|^2+z^2)^{1/2}-z)e^{-\epsilon z}  = |\beta|^2, \\
&& r= \sqrt{|\alpha|^2|\beta|^2+z^2}.
\eea
Then for $\epsilon \geq 0$,
\eqref{eqn:Metric-n=1} determines a hyperk\"ahler metric on $\bC^2$.
For $\epsilon < 0$,
\eqref{eqn:Metric-n=1} determines a hyperk\"ahler metric in the region
defined by $|\alpha \beta| < a$ in $\bC^2$,
and negative of a hyperk\"ahler metric in the region
defined by $|\alpha \beta| > a$ in $\bC^2$.
\end{theorem}

\section{Concluding Remarks}

\label{sec:Conclusions}

Based on the computations in an early paper \cite{Zhou-GH},
we consider the Hessian geometry of toric multi-Taub-NUT metrics.
This means to find the moment maps together with their images,
and complex potential functions on the moment map images,
and use them to introduce some local complex coordinates to express the metrics
and the K\"ahler forms.
Then in the same fashion as in \cite{Zhou-GH},
one can discuss the phase change of the multi-Taub-NUT metrics
by allowing the parameters $c_j$'s to be complex numbers.
Since there is little difference in the treatments,
we have omitted a discussion on this, but instead
focus on the phase transition associated with the parameter $\epsilon$.
Explicit expressions when $\epsilon$ is a nonzero real number
are no longer possible so for explicit examples
we have focused on the case of $n=1$.
In this example
we present a formula for a family of hyperk\"ahler metrics on $\bC^2$ in \eqref{eqn:Metric-n=1}.
The discovery of this formula is a demonstration of the power of Hessian geometry.

We can also take the parameter $\epsilon$ to be a complex number.
This will lead us again to complexified K\"ahler forms as in \cite{Zhou-GH}.
We do not go into this because we plan to write a separate paper to treat
complexifications and phase transitions of K\"ahler forms more generally.

\vspace{.2in}

{\em Acknowledgements}.
The research in this work is partially supported by NSFC grant 11661131005.

\end{document}